\documentclass[a4paper,12pt]{amsart}
  \usepackage{MacrosAnglais}

\begin{document}


\title[Differentiable compactifications]{Symmetric spaces of higher rank do not admit
       differentiable compactifications}
\author{Beno\^{\i}t Kloeckner}

\begin{abstract}
Any nonpositively curved symmetric space 
admits a topological compactification, namely the
Hadamard compactification. For rank $1$ spaces,
this topological compactification can be endowed with
a differentiable structure such that the action of
the isometry group is differentiable. Moreover, the
restriction of the action on the boundary leads to
a flat model for some geometry (conformal, CR or
quaternionic CR depending of the space).
One can ask whether such a differentiable compactification
exists for higher rank spaces, hopefully leading to some knew 
geometry to explore. In this paper we answer negatively.
\end{abstract}

\maketitle

\section{Introduction}

Let $M$ be a symmetric space of nonpositive 
curvature, $G$ its group of isometries and
$G_0$ the identity component
in $G$. 

As a Riemannian manifold, $M$ is a Hadamard space 
and is diffeomorphic to an open ball. Its 
\defini{Hadamard compactification} (or \defini{geodesic
compactification}) is a topological 
gluing of $M$ and its 
\defini{Hadamard boundary} $M(\infty)$ such that
$\adherence{M}=M\cup M(\infty)$ is a closed ball.  
The group $G$ acts continuously on $\adherence{M}$. 

When $M$ is of rank one, that is to say when it is
negatively curved, this topological 
compactification admits ``nice'' models,
carrying an invariant differentiable structure:
in these models, the action of $G$ is differentiable 
on $\adherence{M}$.
Moreover, the restriction of this action to the
boundary is a flat model for some geometry. The 
boundary spheres of the real, complex and quaternionic
hyperbolic spaces wield the standard 
conformal, CR and quaternionic CR structures respectively.
Concerning the octonionic 
hyperbolic plane, the corresponding geometry has not
been studied yet, as far as we know.

It is natural to ask whether such a differentiable
compactification exists when $M$ is of higher rank. One
could expect such a model to give birth to a new,
luckily interesting, geometry.

In this paper, we give a negative answer to this question,
and show that the obstruction comes from the spherical
building at infinity. This combinatorial structure is trivial
only in the Euclidean and rank one spaces. Thus there is
an alternative: a symmetric space of nonpositive curvature
admits either an interesting building at infinity or a 
differentiable compactification, not both.

\subsection{Differentiable compactifications}

Our goal is to extend the differentiable structure of $M$
to the manifold with boundary $\adherence{M}$, 
so that we do not lose symmetry in the process.
This leads to the following definition.

\begin{definition}\label{DefComp}
A \defini{differentiable Hadamard 
compactification} of $M$ is
a differentiable ($\diffb{1}$) 
structure $\ron{D}$ on $\adherence{M}$ compatible with the
differentiable structure of $M$ and such that
the action of $G$ is \diffb{1}.

We can define a \defini{\diffb{r}
Hadamard compactification} in the same way, where
$r$ can be finite, $\infty$ or $\omega$, 
$\diffb{\omega}$ meaning real analytic.
When no precision is given, differentiable
means \diffb{1}.
\end{definition}

By a $\diffb{r}$ action, we mean that
the map $G\times M\to M$ is $\diffb{r}$. It implies
that $G$ acts by $\diffb{r}$ diffeomorphisms and that the
map $G\to\diff{M}$ is continuous in the $\diffb{r}$ topology.

This condition can be greatly relaxed thanks to
the Bochner and Montgomery theorem: if $G$ acts continuously
by $\diffb{r}$ diffeomorphisms, then its action is in fact
$\diffb{r}$ \cite{Bochner-Montgomery}.

For the sake of brevity we will often write ``differentiable compactification''
instead of ``differentiable Hadamard compactification''. However, a symmetric 
space admits other topological
compactifications than the Hada\-mard one, e.g. Martin, Satake and
Furstenberg compactifications. 
It would be interesting to extend our study to these, but the Hadamard 
compactification seems to be of utmost importance for our question.
First, it is very natural, defined directly by the geometry of the
space for a large class of Riemannian manifolds. Second,
there is as far as we know little hope to get a manifold with boundary from
the other compactifications. Either the infinity does not have the right 
dimension (e.g. the Poisson boundary)
or the most natural differentiable structure is that of a manifold with
boundary \emph{and corners} (e.g. the maximal Satake compactification). 
A detailed account on all classical compactifications can be found in
\cite{Borel-Ji} and \cite{Guivarch}.

\subsection{Existence of differentiable compactifications}

Let us now discuss the existence of differentiable compactifications for
the three types of nonpositively curved symmetric spaces.

\paragraph{Symmetric spaces of rank 1.}

It is well known that the real hyperbolic space
$\mH^n$ admits a differentiable Hadamard compactification, given for
example by the closure of Klein's ball:
the central projection of the hyperboloid $Q=-1$ 
(where
$Q$ is the canonical Lorentzian metric on 
$\mR^{n+1}$) gives an embedding of $\mH^n$ into
$\mR\mP^n$ where the group $\SOopq{1}{n}$ of
isometries of $\mH^n$ acts analytically. This 
construction can be generalized
to all symmetric spaces of nonpositive curvature 
and rank $1$.

It is worth noticing that $\mH^n$ admits other differentiable Hadamard
compactifications. For example, the action of $\SOopq{1}{n}$ on
Poincar\'e's ball extends analytically to the closed ball and the resulting
action is not \diffb{1} conjugate to the previous one (this can
be seen by looking at asymptotic geodesics: they are tangent one to another
in the closure
of Poincar\'e's ball, not in Klein's ball.) Details are given in \cite{Kloeckner2},
where it is shown that $\mH^n$ admits an infinite number
of nonconjugate analytic compactifications in the sense of definition
\ref{DefComp}.

\paragraph{Euclidean spaces.}

If $M$ is a Euclidean space, once again it admits a
differentiable Hadamard compactification we briefly
describe.
Identify $\mR^n$ with the affine hyperplane $\{x_0=1\}$ 
of $\mR^{n+1}$ where $n$ is the dimension of $M$.
The projection of center $0$ of $M$ on the open 
upper unit half-sphere is a 
diffeomorphism.
Pushing forward by this map we get an action of 
$G$ (the affine group) on the open upper half-sphere 
whose continuous prolongation to the closed
half-sphere is real analytic. This action is
a real analytic Hadamard compactification of 
$M=\mR^n$.

\paragraph{Symmetric spaces of higher rank.}

The main result of this paper is the following.

\begin{theorem}\label{theoreme}
No noneuclidean symmetric space of rank 
$k\geqslant 2$ admits a differentiable Hadamard
compactification.
\end{theorem}

\paragraph{Structure of the paper.}
From now on, $M$ is supposed to be a symmetric space of rank 
$k\geqslant 2$.

We shall start with a simple remark about the 
natural projection of a fiber $S_xM$ of the unit 
tangent bundle of $M$ on $M(\infty)$. 

In the second section we prove that $\mH^2\times\mR$
admits no differentiable Ha\-da\-mard 
compactification.

Next we generalize this fact to every product 
$F\times\mR^{k-1}$ where $k\geqslant 2$ and
$F$ is a symmetric space of noncompact type of rank
$1$.

Finally, we prove Theorem \ref{theoreme}.

Note that the different parts are more or less 
independent: the proof of Theorem \ref{theoreme} does not
make use of preceding results. However some arguments of 
Section \ref{FxRn} will be useful, and Section \ref{HxR}
gives a good insight of the general phenomenon
on the simplest case.

\section{Apartments and the visual projection}

\subsection{Apartments}

We give some basic vocabulary about the building
structure of $M(\infty)$. More details
can be found in \cite{BallmannGromov}, Appendix 5.
Our main reference for the building structure of
a symmetric space is \cite{Eberlein}.
For details about general buildings, see \cite{Brown}.

Let $A$ be a maximal flat (\lat{i.e.} a totally 
geodesic submanifold isometric to the Euclidean space
of maximal dimension) of $M$, $\adherence{A}$ its closure 
in $\adherence{M}$ and 
$A(\infty)=\adherence{A}\cap M(\infty)$ its
boundary. $A(\infty)$ is called an 
\defini{apartment} of $M(\infty)$. It is a topological 
submanifold.

Every point of $M(\infty)$ belongs to at least
one apartment. A point is 
said to be \defini{regular} if it belongs to exactly 
one apartment, otherwise it is said to be 
\defini{singular}. 

Let $x$ be a point of $M(\infty)$. We denote
by $a(x)$ the set of all apartments containing
$x$. If $x$ is singular, it is said
to have \defini{index} $1$ if $a(x)$ is minimal
with respect to inclusion among sets $a(y)$ of
singular $y$'s.

The connected component of $x$ in the set of 
points $y$ such that $a(x)=a(y)$ is a 
\defini{facet}. Facets are topological submanifolds.
If $x$ is regular, we call its facet 
a \defini{Weyl chamber} or simply a \defini{chamber}; 
if $x$ is singular of index 
$1$, we call its facet a \defini{panel}. 

The dimension of every apartment is $k-1$ (where $k$ is
the rank of $M$). Chambers have
dimension $k-1$, panels have dimension $k-2$.

Two facets are \defini{adjacent} if their
closures intersect. If they are adjacent and of different
dimensions, one is contained in the closure of the other.

The facets form a simplicial complex on
$M(\infty)$ if $M$ is of noncompact type. If
$M$ has a Euclidean factor, some of the cells
are spheres rather than simplicies.

This complex has the incidence structure of a 
spherical thick building, which means:
\begin{enumerate}
\item each apartment is a spherical Coxeter complex (see
      \cite{Brown} for details),
\item for any two facets, there is an apartment containing
      both of them,
\item there exists at least three chambers adjacent to any
      given panel,
\item if there are two apartments $A$, $A'$ containing two
      facets $F$ and $F'$, then there is an isomorphism 
      $A\mapsto A'$ fixing $F$ and $F'$ pointwise.
\end{enumerate}

The group $G$ acts by isomorphisms on this building: it
preserves the adjacency relation and sends facets onto
facets of the same dimension.

\subsection{Non smoothness of the visual 
projection}\label{VP}

Let $x$ be a point of $M$. A unit vector
$v$ tangent to $M$ at $x$ defines a geodesic ray 
$\gamma_v$, hence a point $\gamma_v(\infty)$ of 
the Hadamard boundary $M(\infty)$.
The map 
\[
\pi_x: \begin{array}{rcl}
 S_xM &\to     & M(\infty) \\
   v  &\mapsto & \gamma_v(\infty) \end{array}
\]
is called the \defini{visual projection} from the point $x$.

For all $x$, the visual projection from $x$ is
a homeomorphism. It seems reasonable to expect the
visual projections to be diffeomorphisms for a 
``good'' differentiable Hadamard compactification.
However, it cannot be.

\begin{proposition}\label{RemarqueSimple}
If $M$ is a nonpositively curved symmetric space
of higher rank,
there is no differentiable structure
on $M(\infty)$ such that
all apartments are submanifolds.
\end{proposition}

\begin{proof}
Let $\gamma$ be some geodesic ray in $M$
that is singular of index $1$. 
As the spherical building at infinity
of $M$ is thick, there are at least three (in fact,
an infinite number of) chambers
$C_1$, $C_2$, $C_3$ adjacent to $P$.

For each pair $C_i,C_j$ ($i\neq j$) there is a
flat $A_{ij}$ such that $C_i\subseteq 
A_{ij}(\infty)$ and $C_j\subseteq A_{ij}(\infty)$.
But $A_{ij}(\infty)$ is an embedded 
submanifold of $M(\infty)$, thus $C_i$ and $C_j$
have opposite tangent half spaces $E_i$, $E_j$ at 
$\gamma(\infty)$. See figure \ref{RS1}.

\begin{figure}[htb]\begin{center}
  \includegraphics[width=90mm]{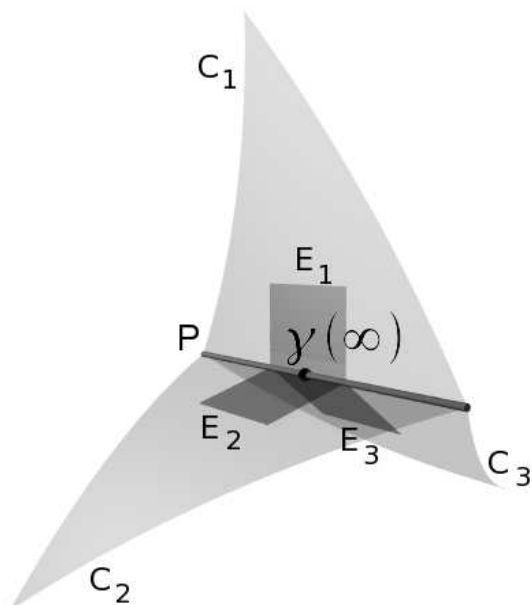}
  \caption{Three chambers meeting at a panel.}
  \label{RS1}
\end{center}\end{figure}

Thus we get three half subspaces $E_1$, $E_2$, $E_3$ 
of $T_{\gamma(\infty)}M(\infty)$ such that
$E_1=-E_2$, $E_1=-E_3$, $E_2=-E_3$, a 
contradiction.
\end{proof}

\begin{corollary}
There is no differentiable structure 
on $M(\infty)$ such that
$\pi_x$ is a diffeomorphism for all $x\in M$.
\end{corollary}

\begin{proof}
Suppose there is such a differentiable 
structure.

Let $A$ be a maximal flat of $M$, $x$ be a point 
of $A$.

Then $S_xA$ is an embedded submanifold of $S_xM$
and $\pi_x$ is a diffeomorphism. Thus 
$A(\infty)=\pi_x(S_xA)$ is a submanifold of 
$M(\infty)$. A contradiction with
Proposition \ref{RemarqueSimple}.
\end{proof}

\section{Study of $\mH^2\times\mR$}\label{HxR}

We summarize briefly the building structure 
of $\mH^2\times\mR$. 

\begin{figure}[hbt]\begin{center}
  \includegraphics[width=90mm]{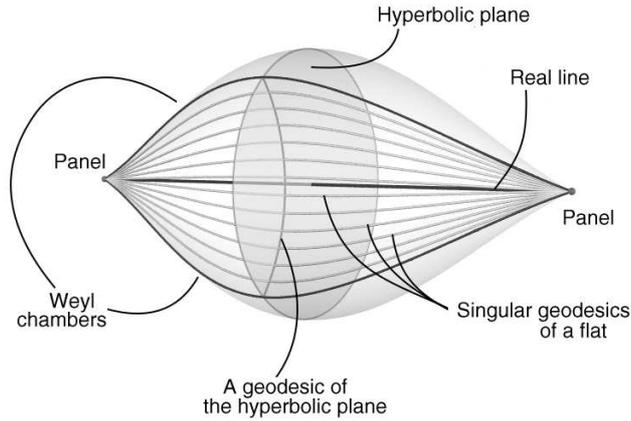}
  \caption{A family of parallel singular geodesics lying on the
           same flat.}
  \label{RS2}
\end{center}\end{figure}

The singular geodesics are those of the form
$\{x\}\times\mR$ where $x$ is a point of $\mH^2$; they are
parallel (asymptotic at both ends) to one another.
The maximal flats are the products
$\gamma\times\mR$ where $\gamma$ is a geodesic of $\mH^2$
(see figure \ref{RS2}).

\begin{figure}[htb]\begin{center}
  \includegraphics[width=90mm]{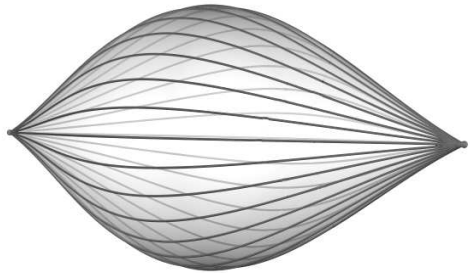}
  \caption{Here we show the Weyl chambers in the boundary of
           $\mH^2\times\mR$. Any choice of two of them corresponds
           to a geodesic of $\mH^2$, a flat of $\mH^2\times\mR$
           and an apartment.}
  \label{RS3}
\end{center}\end{figure}

The boundary of $\mH^2\times\mR$ is a $2$-sphere partitionned into
two points and a family of nonintersecting curves joining them. 
The points are the end points of every singular geodesic, therefore
panels of the building. The curves are the Weyl chambers,
there is one of them for each point in the boundary
of $\mH^2$ (see figure \ref{RS3}). The union of any two of them and of the
two panels is an apartment.

\begin{proposition}\label{HDR}
The space $\mH^2\times\mR$ admits no 
differentiable
Hadamard compactification.
\end{proposition}

\begin{proof}
Suppose $\ron{D}$ is a differentiable Hadamard
compactification of $M=\mH^2\times\mR$ and 
$\adherence{M}$ is endowed with $\ron{D}$.

Let $\gamma=\{x\}\times\mR^+$ be any singular geodesic 
ray of unit speed.
We denote by $\gamma(\infty)$ the point of
$M(\infty)$ defined by $\gamma$.

Since $\gamma(\infty)$ is fixed by all 
orientation-preserving isometries of
$\mH^2$, the deri\-va\-tives
of these isometries give a linear 
representation of $\PSLDR$ on 
$T_{\gamma(\infty)}\adherence{M}$.
This representation is reducible as
$T_{\gamma(\infty)} M(\infty)$ is an invariant
subspace. Let $\rho$ be the induced representation.
As a representation of a simple Lie group 
$\rho$ is trivial or faithful.

Let $s_x$ be the geodesic symmetry of $\mH^2$ 
around $x$. We identify $s_x$ with the isometry
$s_x \times \id$ of $\mH^2\times\mR$.
For every time $t\in\mR$, $ds_x(\gamma(t))$ has
eigenvalues $1$, $-1$, $-1$. By continuity,
$ds_x(\gamma(\infty))$ must have the same 
eigenvalues, thus the restriction to 
$T_{\gamma(\infty)} M(\infty)$ of 
$ds_x(\gamma(\infty))$ is $-\id$, that
is to say $\rho(s_x)=-\id$.

Let $y$ be any point of $\mH^2$ different from $x$.
Then $\rho(s_y)=-\id$ too. Thus $\rho(s_xs_y)=\id$,
but $s_xs_y$ is a non-trivial hyperbolic 
transformation (it is a translation along the 
geodesic containing $x$ and $y$). Thus $\rho$
is neither faithful nor trivial, a contradiction.
\end{proof}

\section{Product of a Euclidean space by a rank 
         $1$ space}\label{FxRn}

We now generalize Proposition \ref{HDR} to the case
when $M=F\times\mR^{k-1}$ is the product of an
Euclidean space by a symmetric space $F$ 
of rank $1$. 

However, as we will need it later, we prove 
something stronger.

\begin{definition}
A \defini{weak differentiable 
Hadamard compactification} of $M$ is defined 
as a differentiable Hadamard compactification where
we replace $G$ by its identity component $G_0$. 
\end{definition}

We prove that $M$ admits no weak differentiable
Hadamard compactification. Thus, in order to
generalize the argument used in the proof of
Proposition \ref{HDR} we need the geodesic
symmetries to belong to $G_0$. Of course this is false
if $F$ is of odd dimension, and we shall
use another argument in this case.

\begin{proposition}
Let $F$ be a rank $1$ symmetric space of noncompact 
type. If $\dim{F}$ is even, then the geodesic
symmetries belong to $G_0$. If $\dim{F}$ is
odd, then $F$ is a real hyperbolic space 
$\mH^{2m+1}$.
\end{proposition}

\begin{proof}
From the classification of symmetric spaces (see 
for example \cite{Helgason}, Chapter IX) we know that the 
rank $1$
symmetric spaces of noncompact type are: the
real hyperbolic spaces, the complex hyperbolic
spaces, the quaternionic hyperbolic spaces and
an exceptionnal space, the octonionic hyperbolic 
plane (therefore the last assertion is clear). The 
identity components of their isometry groups are
respectively 
$\SOopq{1}{n}$, $\SUpq{1}{n}$, $\Sppq{1}{n}$
and $F_{4(-20)}$, which are simple Lie groups. 

We shall use the following criterion.
Let $\al{g}=\al{k}+\al{p}$ be a Cartan decomposition
of the Lie algebra of $G_0$. Then
the geodesic symmetries are in $G_0$ if and only 
if $\al{k}$ contains a maximal abelian 
algebra of $\al{g}$ (see \cite{Helgason}, Chapter IX �3).

It is know sufficient to compare the ranks of
$\al{g}$ and $\al{k}$ (most of them can be found
in \cite{Knapp2}, Appendix C):
\begin{itemize}
\item in the real case, $\al{g}=\sopq{n}{1}$ is of rank
      $\lfloor\frac{1}{2}(n+1)\rfloor$ and
      $\al{k}=\so{p}$, of rank $\lfloor n/2\rfloor$.
      These two ranks coincide exactly when $n$ is even,
\item in the complex case, $\al{g}=\supq{n}{1}$ and
      $\al{k}=\al{s}(\up{n}\times\up{1})$ are both of rank
      $n$,
\item in the quaternionic case, $\al{g}=\sppq{n}{1}$ and
      $\al{k}=\spp{n}\times\spp{1}$ are both of rank $n+1$,
\item in the octonionic case, $\al{g}=\al{f}_{4(-20)}$ and
      $\al{k}=\al{so}(9)$ are both of rank $4$. I wish to thank
      Fokko du Cloux, J\'{e}r\^{o}me Germoni and Bruno S\'{e}vennec
      for explaning this case to me.
\end{itemize}
\end{proof}

We can now prove the following.

\begin{proposition}\label{produit}
If $M=F\times\mR^{k-1}$ where $F$ is a
rank $1$ symmetric space of noncompact type,
then $M$
admits no weak differentiable Hadamard compactification.
In particular, $M$ admits no differentiable Hadamard
compactification.
\end{proposition}

\begin{proof}
Suppose there is such a differentiable structure
on $\adherence{M}$.

We identify an isometry $g$ of $F$ with the isometry
$g\times \id$ of $M$. We denote the component of 
identity of the group of isometries of $F$  by 
$G_0^F$ and consider it a subgroup of $G_0$. 

We first suppose that $F$ is of even dimension.
 
Let $s_x$ and $s_y$
be the geodesic symmetries around two different
points $x$ and $y$ of $F$. Let $d$ be a geodesic
of $\mR^{k-1}$. Then $\gamma_1=\{x\}\times d$ and
$\gamma_2=\{y\}\times d$ are asymptotic and
$z=\gamma_1(+\infty)=\gamma_2(+\infty)$ is a singular
point of index $1$ of $M(\infty)$.

Differentiation $g\mapsto dg(z)$ gives us
a linear representation
of $G_0^F$ in $T_{z} \adherence{M}$. This 
representation is reducible
since $T_{z} M(\infty)$ is an invariant 
subspace. Let $\rho$ be the representation induced on
$T_{z} M(\infty)$.

We now decompose the representation $\rho$.
For all $t\in\mR$, the eigenvalues of
$ds_x(\gamma_1(t))$ are $1$ with multiplicity
$k-1$ and $-1$ with multiplicity $\dim{F}$.
Thus, $\rho(s_x)$ must have eigenvalues $1$ with 
multiplicity $k-2$ and $-1$ with multiplicity 
$\dim{F}$ (an eigendirection
transverse to the boundary must have a nonnegative 
eigenvalue).

We shall decompose $\rho$ using the
following lemma.

\begin{lemma}\label{lemm1produit}
The panel $P$ of $\gamma_1(\infty)$ is a 
differentiable submanifold of $M(\infty)$.
\end{lemma}

\begin{proof}
The panel $P$ is pointwise
fixed by all isometries of $F$. Thus it is pointwise 
fixed by $s_x$. In a local chart, it is defined by 
$p\in P\Rightarrow s_x(p)-p=0$. Since $s_x-\id$ has 
rank $\dim{F}$ and $\dim{M(\infty)}=\dim{F}+k-2$, the
inverse function theorem implies that $P$ 
is contained in a differentiable submanifold of 
$M(\infty)$ of dimension 
$k-2$ (namely the set of fixed points of $s_x$). 
But, as a panel in the boundary of a rank $k$ 
symmetric space, it is an open topological manifold 
of dimension $k-2$.
Thus $P$ is a differentiable submanifold of dimension 
$k-2$ of $M(\infty)$.
\end{proof}

The tangent space $T_z P$ is an invariant
subspace of $\rho$. Thus, this representation splits 
in two 
parts : $\rho=\rho_0\oplus\rho_1$ where $\rho_0$ is 
the trivial representation of dimension $k-2$ and 
$\rho_1$ is a representation of dimension $\dim{F}$.

Now we have a representation $\rho_1$ (which, as 
$G_0^F$ is simple, must be faithful or trivial) 
with $\rho_1(s_x)=-\id$. So $\rho_1$ cannot be 
trivial. But $ds_y(\gamma_2(t))$ has the same
eigenvalues as $ds_x(\gamma_1(t))$, and thus 
$\rho_1(s_y)=-\id$ too. Now we have 
$\rho_1(s_xs_y)=\id$ with $s_xs_y$ a hyperbolic
transformation, so $\rho_1$ cannot be faithful, a 
contradiction.

Suppose now that $F$ is of odd dimension. 

We have
$F=\mH^{2m+1}$, the real hyperbolic space. The 
geodesic symmetries are not in $G_0$ (their 
determinant is $-1$) and
we shall use Proposition \ref{RemarqueSimple}.

Let $A$ be a maximal flat. Then $A$ is the product of 
a geodesic $\gamma$ of $F$ by $\mR^{k-1}$. Let 
$r\in G_0^F$ be 
the rotation of angle $\pi$ around $\gamma$ in $F$; 
$A$ is the set of fixed points of $r$, 
and $A(\infty)$ is the set of fixed
points of $r$ in $M(\infty)$.
Since $r$ is an involution, $A(\infty)$ is a
submanifold
of $M(\infty)$, a contradiction to 
Proposition \ref{RemarqueSimple}.

We give an alternative proof for the odd dimension 
case, less elegant but more useful for the proof
of Theorem \ref{theoreme}. We can define the 
representation $\rho_1$ like in the even dimension
case. Then $\rho_1$ is a representation
of dimension $\dim{F}$ of $G_0^F$.
Since, for all $x\in M$ fixed by $r$, $dr(x)$ has 
eigenvalues $-1$ with multiplicity 
$2m$ and $1$ with multiplicity $k$,
$\rho_1(r)$ has eigenvalues $-1$ with multiplicity
$2m$ and $1$ with multiplicity $1$. Thus
$\rho_1(r)\neq\id$, hence $\rho_1$ is not trivial.

But $G_0^F=\SOopq{2m+1}{1}$ admits
no non trivial representation of dimension
less than $2m+2$, a contradiction.
\end{proof}

\section{Proof of Theorem \ref{theoreme}}\label{Stheoreme}

We shall now prove Theorem \ref{theoreme} with the
same ideas that we used for the previous
propositions.

Let $M$ be a noneuclidean symmetric space of 
nonpositive curvature of dimension $n$ and rank $k>1$. 
As before,
$G$ is the group of all isometries of $M$,
$G_0$ is the identity component of $G$
and $\al{g}$ is the Lie algebra of $G$.

Suppose that there exists a differentiable Hadamard
compactification $\ron{D}$ of $M$.

We denote by $\alpha$ the action of $G$ on 
$\adherence{M}$. We also denote by $\alpha$ 
the corresponding action of $\al{g}$.

The first step is to find in $M$ an embedded
product $F\times \mR^{k-1}$.

Let $\gamma$ be a singular geodesic of index $1$.
Let $F_{\gamma}$ be the union of geodesics parallel
to $\gamma$ (recall that parallel means that they are both positively
and negatively asymptotic). Then, $F_{\gamma}$ is a totally
geodesic submanifold of $M$ isometric to a product
$F\times\mR^{k-1}$ where $F$ is a symmetric space of 
rank $1$ (see \cite{Eberlein} Section 2.11).

Let $\adherence{F}_{\gamma}$ be the closure of
$F_{\gamma}$ in $\adherence{M}$ and 
$F_{\gamma}(\infty)=
\adherence{F}_{\gamma}\cap M(\infty)$.

It would be interesting to prove that 
$\adherence{F}_{\gamma}$ is a submanifold of 
$\adherence{M}$, since we could directly use
Proposition \ref{produit} to get a contradiction,
but a weaker statement (namely Lemma \ref{lemm5theorem})
will be sufficient.

Up to a change
of parametrization we can write
\[\gamma(t)=(p,(t,0,\dots,0))\]
where $p\in F$, and $F^t=F\times\{(t,0\dots,0)\}$ is
identified with its embedding into $M$.

Since $F$, identified with $F^0$,
is a totally geodesic submanifold of $M$,
the Lie algebra of the group $G^F$ of isometries of $F$ is
a subalgebra of $\al{g}$ and the identity component $G_0^F$
of $G^F$ is a subgroup of $G_0$ (that's why we needed
the stronger statement in Proposition \ref{produit}). Thus
taking derivatives gives us a 
representation $\rho$ of $G_0^F$ on 
$T_{\gamma(\infty)} \adherence{M}$.

Let 
$\al{k}^t\oplus\al{p}^t$ be the Cartan decomposition
of $\al{g}$ at $\gamma(t)$.
Since $F_{\gamma}$ is a totally geodesic submanifold
we have a further decomposition
$\al{p}^t=\al{p}^t_F\oplus\al{p}^t_{eucl}\oplus\al{p}^t_0$
where the terms are pairwise orthogonal (with respect to
the Killing form),
$\al{p}^t_F$ is mapped by $\alpha$
onto $T_pF^t$, $\al{p}^t_{eucl}$ is mapped onto 
$T_p\mR^{k-1}$ and $\al{p}^t_0$ is mapped
onto $\ortho{(T_pF_{\gamma})}$. We define
$\al{p}^t_{\gamma}=\al{p}^t_F\oplus\al{p}^t_{eucl}$;
$\alpha$ maps $\al{p}^t_{\gamma}$ onto
$T_{\gamma(t)} F_{\gamma}$. Moreover,
$\al{p}^t_{\gamma}$ is the set of all Killing fields in
$\al{p}$ commuting with that Killing field $X\in\al{p}^t$ such that
$\alpha(X)$ is the unit tangent vector of
$\gamma$.

We shall split $\rho$ in three parts in correspondence
with the splitting 
$\al{p}^t=\al{p}^t_F\oplus\al{p}^t_{eucl}\oplus\al{p}^t_0$.
To achieve this, we use the following stability result.

\begin{lemma}\label{lemmrep}
Let $K$ be any compact group and $(\mu_t)_{t\in\mR}$ be a continuous
family of linear representations of $K$ on some finite-dimensional
real vector space
$V$. Then for all pairs $(t_1,t_2)$ of real numbers, the representations
$\mu_{t_1}$ and $\mu_{t_2}$ are conjugate.
\end{lemma}

\begin{proof}
As $K$ is compact and $V$ is finite-dimensional, the conjugacy 
class of a representation $\mu_t$ is determined by its character. 
More precisely, the multiplicity in $\mu_t$ of some irreducible 
representation $\nu$ is given by the scalar product of the 
characters of $\mu_t$ and $\nu$, hence is a continuous map. This 
multiplicity is an integer and is thus constant.
\end{proof}

\begin{lemma}\label{lemm1theorem}
The tangent subspaces $T_{\gamma(t)} F^t$,
$T_{\gamma(t)} \mR^{k-1}$ and 
$\ortho{(T_{\gamma(t)} F_{\gamma})}$ admit limits
when $t\rightarrow\infty$, denoted respectively
by $V_F$, $V_{eucl}$ and $V_0$.
Moreover one has
$T_{\gamma(\infty)}\adherence{M}=V_F\oplus V_{eucl}\oplus V_0$.
\end{lemma}

\begin{proof}
Let $K_0^F$ be the isotropy group of $p$ in $F$.
Taking derivatives gives linear representations
$\rho^t$ of $K_0^F$ in $T_{\gamma(t)} M$ for all $t\in\mR$.
Then $T_{\gamma(t)} F^t$, $T_{\gamma(t)}\mR^{k-1}$ and
$\ortho{(T_{\gamma(t)}F_\gamma)}$ are invariant spaces
of $\rho_t$. 

By continuity, the restriction $\rho^{\infty}$ of $\rho$ to
$K_0^F$ splits into three parts and
the conclusion holds.
\end{proof}

From Lemma \ref{lemmrep} we deduce that the action of
$\rho^{\infty}$ is conjugate with that of $\rho^0$.
Since $\rho^{t}$ acts trivially on 
$T_{\gamma(t)}\mR^{k-1}$, $\rho^{\infty}$ acts trivially 
on $V_{eucl}$.

We shall now prove that $V_F$ is an invariant subspace for
$\rho$.

\begin{lemma}\label{lemm2theorem}
Let $\ron{O}$ be the orbit of $\gamma(\infty)$ under the
action of $G_0$. Then 
$V_0=T_{\gamma(\infty)} \ron{O}$ and, for all $t$, the restriction
of $\alpha_{\gamma(\infty)}$ to $\al{p}^t_0$ 
is one-to-one and onto $V_0$.
\end{lemma}

\begin{proof}
As an orbit, $\ron{O}$ is a submanifold of $M(\infty)$,
invariant under the action of $G_0$, thus 
$T_{\gamma(\infty)}\ron{O}$ is an invariant space
of $\rho$. 

The definition of
$\ron{O}$ shows that $\alpha$ sends $\al{p}^t_0$ on 
$T_{\gamma(\infty)}\ron{O}$. We want to prove that
this map is one-to-one and onto.

Let $H$ be an element of $\al{p}^t_0$. By definition,
$\alpha(\exp(H))(\gamma)$ is not parallel to $\gamma$. If 
$\alpha(H)_{\gamma(\infty)}=0$, then
$\alpha(H)_{\gamma(-\infty)}\neq 0$. But after
conjugacy by the geodesic symmetry at $\gamma(t)$ we find
$\alpha(-H)_{\gamma(-\infty)}=0$, a contradiction. Thus
the restriction of $\alpha_{\gamma(\infty)}$ to $\al{p}^t_0$
is one-to-one.

Since
$\dim{\al{p}^t_0}=n-\dim{F_{\gamma}}=
(n-1)-(\dim{F_{\gamma}}-1)=\dim{\ron{O}}$, it is onto
and $V_0=T_{\gamma(\infty)}\ron{O}$.
\end{proof}

\begin{lemma}\label{lemm3theorem}
The subspace $V_F$ contains no subspace where $\rho$
is trivial.

Let $X$ be that vector of $\al{p}^0_{eucl}$
such that $\gamma'(t) = \alpha(X)_{\gamma(t)}$. 
The linear operator 
$d_{\gamma(\infty)} \alpha(\exp X)-\id$ acting on
$V_0$ is of maximal rank.

\end{lemma}

\begin{proof}
The linear action $\rho^{\infty}$
 of $K_0^F$ on $T_{\gamma(\infty)} M(\infty)$ is
transitive on $V_F$ and $\rho^{\infty}$ is the restriction
of $\rho$ to $K_0^F$ thus $V_F$ contains no trivial part.

To prove the second part of the lemma, we use the root space decomposition
$\al{g}=\al{g}_0+\sum \al{g}_{\lambda}$ given by
some maximal flat containing $\gamma$. 
We have 
\[\al{p}^0_0\subseteq \sum_{\lambda(X)\neq 0} \al{g}_{\lambda}\]
and $Ad(\exp X) = e^{\lambda(X)} I$ on $\al{g}_{\lambda}$.

Since $\alpha$ is onto from $\al{p}^0_0$ to $V_0$, it is
onto from $\sum_{\lambda(X)\neq 0} \al{g}_{\lambda}$ to
$V_0$.

For all $H\in \al{g}_{\lambda}$, we have
\[d_{\gamma(\infty)} \alpha(\exp X)(\alpha(H))
  = \alpha(Ad(\exp(X))(H))
  = e^{\lambda(X)}\alpha(H),\]
and thus $d_{\gamma(\infty)} \alpha(\exp X)-\id$ is
nondegenerate on $V_0$.
\end{proof}

\begin{lemma}\label{lemm4theorem}
The panel $P$ of $\gamma(\infty)$ is a submanifold
of $M(\infty)$ and its tangent space at $\gamma(\infty)$
equals $V_{eucl}\cap T_{\gamma(\infty)}M(\infty)$.
\end{lemma}

\begin{proof}
The panel $P$ is contained
in the set of the points of $M(\infty)$ left
fixed by the actions of $G$ and of $\exp(X)$.
Written in local coordinates, this gives us an infinite
system of equations. By Lemma \ref{lemm3theorem}
we know that this system is of rank at least
$d=\dim V_0+\dim V_F$ at $\gamma(\infty)$.
We can extract a subsystem of $d$ equations that is
of maximal rank at $\gamma(\infty)$. The inverse function
theorem implies that this subsystem defines a submanifold
of $M(\infty)$ of dimension 
$k-2=\dim V_{eucl}\cap T_{\gamma(\infty)}M(\infty)$ and 
containing
$P$. But $P$ is topologically a manifold of dimension $k-2$ and
thus must be a differentiable submanifold of $M(\infty)$.

Since $K_0^F$ acts trivially on $P$, its tangent space
must be $V_{eucl}\cap T_{\gamma(\infty)}M(\infty)$.
\end{proof}

Since $G$ acts trivially on $P$, it must preserve
its tangent space. We are now ready to prove the following.

\begin{lemma}\label{lemm5theorem}
The subspace $V_F$ is invariant by $\rho$.
\end{lemma}

\begin{proof}
From previous lemmas we know that $V_0$ and 
$V_{eucl}\cap T_{\gamma(\infty)}M(\infty)$
are invariant subspaces for $\rho$.
Since $\rho$ is totally reducible, there exists some
subspace $V'$ invariant by $\rho$ such that one has
the following decomposition
$T_{\gamma(\infty)} M(\infty) = 
V'\oplus V_{eucl}\cap T_{\gamma(\infty)}M(\infty) \oplus V_0$. 
But $V'$ must be invariant by 
$\rho^{\infty}$, and thus $V'=V_F$ and $V_F$ is invariant by
$\rho$.
\end{proof}

Denote by $\rho_1$ the representation of $G_0^F$ induced
by $\rho$ on $V_F$.

Since the restriction of $\rho_1$ to $K_0^F$ is the limit
of the restriction of $\rho^t$ to $T_{\gamma(t)} F^t$,
we can now use the same arguments as in the proof 
of Proposition \ref{produit}.

If $F$ is of even dimension, $\rho_1$ is neither trivial nor faithful, a contradiction.

If $F$ is of odd dimension, $\rho_1$ is a nontrivial
representation of dimension $\dim{F}$ of $G_0^F$,
a contradiction.

Theorem \ref{theoreme} is proved.
Note that we actually get something stronger: there exists no \insiste{weak}
differentiable Hadamard compactification of $M$; the obstructions to
differentiability appear in the identity component of $G$.

\nocite{Serre, Vinberg, Knapp, Guivarch}
\bibliographystyle{plain}
\bibliography{biblio.bib}

\def\dbar{\leavevmode\hbox to 0pt{\hskip.2ex \accent"16\hss}d}
\begin{thebibliography}{10}

\bibitem{BallmannGromov}
Werner Ballmann, Mikhael Gromov, and Viktor Schroeder.
\newblock {\em Manifolds of nonpositive curvature}, volume~61 of {\em Progress
  in Mathematics}.
\newblock Birkh\"auser Boston Inc., Boston, MA, 1985.

\bibitem{Bochner-Montgomery}
Salomon Bochner and Deane Montgomery.
\newblock Groups of differentiable and real or complex analytic
  transformations.
\newblock {\em Ann. of Math. (2)}, 46:685--694, 1945.

\bibitem{Borel-Ji}
Armand Borel and Lizhen Ji.
\newblock {\em Compactifications of symmetric and locally symmetric spaces}.
\newblock Mathematics: Theory \& Applications. Birkh\"auser Boston Inc.,
  Boston, MA, 2006.

\bibitem{Brown}
Kenneth~S. Brown.
\newblock {\em Buildings}.
\newblock Springer-Verlag, New York, 1989.

\bibitem{Eberlein}
Patrick~B. Eberlein.
\newblock {\em Geometry of nonpositively curved manifolds}.
\newblock Chicago Lectures in Mathematics. University of Chicago Press,
  Chicago, IL, 1996.

\bibitem{Guivarch}
Yves Guivarc'h, Lizhen Ji, and J.~C. Taylor.
\newblock {\em Compactifications of symmetric spaces}, volume 156 of {\em
  Progress in Mathematics}.
\newblock Birkh\"auser Boston Inc., Boston, MA, 1998.

\bibitem{Helgason}
Sigur{\dbar}ur Helgason.
\newblock {\em Differential geometry and symmetric spaces}.
\newblock Pure and Applied Mathematics, Vol. XII. Academic Press, New York,
  1962.

\bibitem{Kloeckner2}
Beno{\^{\i}}t Kloeckner.
\newblock On differentiable compactifications of the hyperbolic space.
\newblock {\em Transform. Groups}, 11(2):185--194, 2006.

\bibitem{Knapp2}
Anthony~W. Knapp.
\newblock {\em Representation theory of semisimple groups}, volume~36 of {\em
  Princeton Mathematical Series}.
\newblock Princeton University Press, Princeton, NJ, 1986.

\bibitem{Knapp}
Anthony~W. Knapp.
\newblock {\em Lie groups beyond an introduction}, volume 140 of {\em Progress
  in Mathematics}.
\newblock Birkh\"auser Boston Inc., Boston, MA, second edition, 2002.

\bibitem{Serre}
Jean-Pierre Serre.
\newblock {\em Alg\`ebres de {L}ie semi-simples complexes}.
\newblock W. A. Benjamin, inc., New York-Amsterdam, 1966.

\bibitem{Vinberg}
{\`E}.~B. Vinberg, editor.
\newblock {\em Lie groups and {L}ie algebras, {III}}, volume~41 of {\em
  Encyclopaedia of Mathematical Sciences}.
\newblock Springer-Verlag, Berlin, 1994.

\end{thebibliography}

\end{document}